\documentclass[sts,preprint]{imsart}

\usepackage{amsthm, amsmath, amssymb, natbib, relsize, mathrsfs}
\RequirePackage[colorlinks,citecolor=blue,urlcolor=blue]{hyperref}
\usepackage{graphicx, subcaption}
\graphicspath{{figures/}} 

\startlocaldefs
\DeclareMathOperator*{\argmax}{arg\,max}
\newtheorem{thm}{Theorem}
\newtheorem*{thm*}{Theorem}

\newtheorem*{corly*}{Corollary}
\newtheorem{lem}{Lemma}[thm]
\newtheorem*{lem*}{Lemma}
\endlocaldefs

\allowdisplaybreaks

\begin{document}

\begin{frontmatter}

\title{Reconciling the Bayes Factor and Likelihood Ratio for Two Non-Nested Model Selection Problems}
\runtitle{Bayes Factor and Likelihood Ratio}

\begin{aug}
  \author{\fnms{Danica M.}  \snm{Ommen}\thanksref{t1,t2}\ead[label=e1]{dmommen@iastate.edu}}
  \and
  \author{\fnms{Christopher P.}  \snm{Saunders}\thanksref{t1}%
  \ead[label=e3]{christopher.saunders@sdstate.edu}%
}

  \thankstext{t1}{Research supported by the National Institute of Justice, Office of Justice Programs, US Department of Justice under Award No. 2014-IJ-CX-K088. The opinions and conclusions or recommendations expressed in this presentation are those of the author and do not necessarily represent those of the Department of Justice.}
  \thankstext{t2}{Research partially supported by the Center for Statistics and Applications in Forensic Evidence (CSAFE) through Cooperative Agreement \#70NANB15H176 between NIST and Iowa State University, which includes activities carried out at Carnegie Mellon University, University of California Irvine, and University of Virginia.}
  \runauthor{Ommen \& Saunders}

  \affiliation{Iowa State University and South Dakota State University}

  \address{Iowa State University, Department of Statistics, Ames, IA 
          \printead{e1};}

  \address{South Dakota State University, Department of Mathematics \& Statistics, Brookings, SD
          \printead{e3}.}

\end{aug}

\begin{abstract}
In statistics, there are a variety of methods for performing model selection that all stem from slightly different paradigms of statistical inference.  The reasons for choosing one particular method over another seem to be based entirely on philosophical preferences.  In the case of non-nested model selection, two of the prevailing techniques are the Bayes Factor and the likelihood ratio. This article focuses on reconciling the likelihood ratio and the Bayes Factor for comparing a pair of non-nested models under two different problem frameworks typical in forensic science, the common-source and the specific-source problem.  We show that the Bayes Factor can be expressed as the expected value of the corresponding likelihood ratio function with respect to the posterior distribution for the parameters given the entire set of data where the set(s) of unknown-source observations has been generated according to the second model. This expression leads to a number of useful theoretic and practical results relating the two statistical approaches.  This relationship is quite meaningful in many scientific applications where there is a general confusion between the various statistical methods, and particularly in the case of forensic science.
\end{abstract}


\end{frontmatter}

\section{Introduction}\label{intro-section}

One of the difficulties with having so many different statistical paradigms is trying to choose which method best suits your problem.  In many cases, there are philosophical arguments that support more than one method, so the scientist will need to choose.  During the process of developing and refining methods for a particular application, some scientists may choose one method while others choose an opposing method.  In the case of non-nested model selection, both the likelihood ratio, most often associated with the classical paradigm, and the Bayes Factor, most often associated with the Bayesian paradigm, are appropriate analyses aimed at answering the same basic question: which model is better supported by the data. Inevitably, the next scientist comes along and must make a choice about which methodology to employ.  Unfortunately, there is often no clear way of directly comparing the quality and rigor of two final products created under different paradigms.  Further complicating matters, we often use ambiguous language and common terms between the two different paradigms to describe the methods.  For example, Bayes Factors and likelihood ratios are both commonly called ``likelihood ratio" when used in forensic science applications, regardless of the statistical paradigm used in the analysis.  This makes it very difficult to compare the performance of one ``likelihood ratio" to another in forensics since the two may not be computed in a similar fashion while both are purportedly expressing the same value. The goal of this article is to take the first steps towards providing subscribers of one model selection paradigm a direct link to the other for two very useful non-nested model selection frameworks.

To begin, Section~\ref{nnms-section} will summarize the two different model selections frameworks that we call the common-source and the specific-source problems.  Then for each of the problems, the two models from which the selection is to be made will be defined.  In Section~\ref{methods}, the forms of both the Bayes Factor and the likelihood ratio for the two non-nested model selection problems will be given.  Finally, the relationships between them are explored in Section~\ref{LRvsBF1}.  The interested reader is directed to the supplementary material \citet{BFvsLRsupplement} for further details that have been omitted from this article for the sake of conciseness and clarity of presentation.

\section{Non-Nested Model Selection}\label{nnms-section}

One of the most common, and increasingly the most difficult, areas of statistical application to forensic science is in the subject of forensic identification of source problems.  The general idea of identification of source is that you have an object related to the perpetration of a crime, and you wish to determine where that object originated.  For example, a fingerprint is left at the scene of a murder, and you want to determine if the print originates from a finger of the suspect.  Similarly, suppose you have two different crime scenes with fingerprints recovered, and you want to know whether the prints were left by the same finger of an unknown perpetrator.  (See Ommen and Saunders \citet{CSvsSS-LPR} for further details of forensic identification of source problems of these types.)  As it turns out, these two scenarios can be expressed statistically as two different non-nested model selection problems.

\subsection{Common-Source Models}\label{cs}

\begin{figure} 
\includegraphics[width=\textwidth]{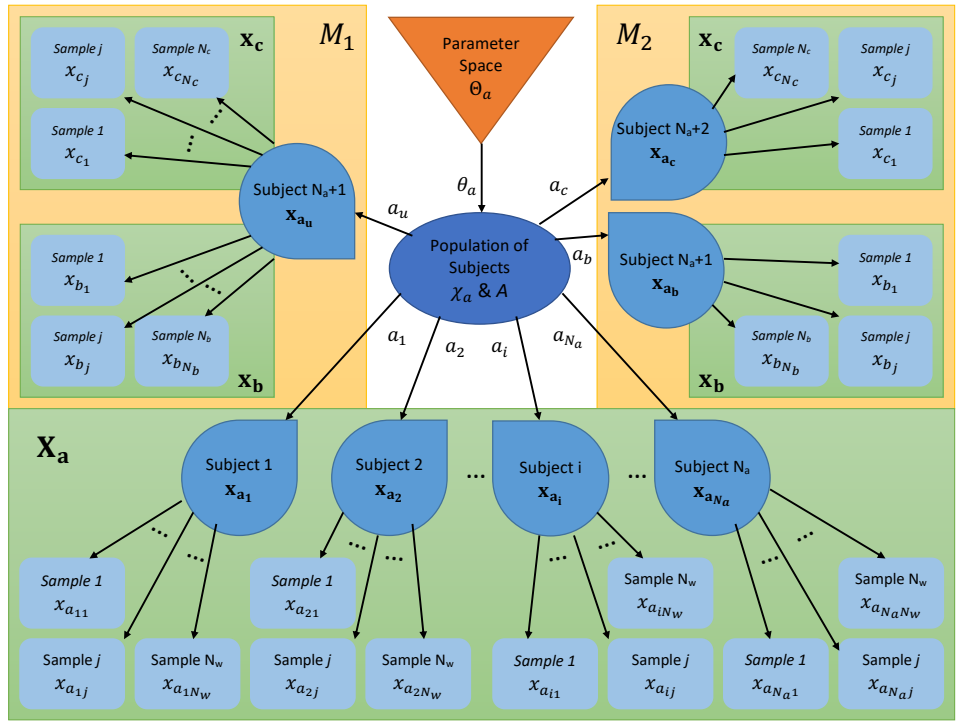}
\caption[]{Hierarchy and notation for the common-source model selection problem where the green boxes denote the datasets and the yellow boxes denote the two competing models.}
\label{CSmodels}
\end{figure}

The idea of the common-source non-nested model selection problem is that you have built up a dataset, $\mathbf{X_a}$, of observations from many different subjects (the number of subjects is denoted $N_a$, and the number of samples from within each subject is denoted $N_w$) in some population and two sets of observations with unknown sources, $\mathbf{x_b}$ and $\mathbf{x_c}$, have been observed (the sample sizes are denoted by $N_b$ and $N_c$, respectively).  We assume that each set of observations has been generated by one single source, and that they are generated from the model which produced the dataset $\mathbf{X_a}$.  Now, it is of interest to determine whether the two sets of unknown-source observations have been generated from a common unspecified (random) subject in the population, denoted by $M_1$, or if the two sets of unknown-source observations have been generated from two different unspecified (random) subjects in the population, denoted by $M_2$.  In the absence of simplifying assumptions regarding the subjects in the dataset $\mathbf{X_a}$, the model describing the generation of $\mathbf{X_a}$ implicitly defines a latent random variable related to the selection of a subject from the population.  This latent random variable will be denoted by $A$ and the observed value of $A$ corresponding to the $i^{th}$ subject in the population will be denoted by $a_i$.  The corresponding sampling models for the entire set of the data, denoted $\mathbf{X_N} = \{ \mathbf{X_a},  \mathbf{x_b}, \mathbf{x_c}\}$ where $N$ denotes the index for the various sample sizes, are given in the supplementary material \citet{BFvsLRsupplement} and visualized in Figure~\ref{CSmodels}.  An example demonstrating the set-up for an application of the common-source framework to forensic evidence is given in Ommen and Saunders \citet{CSvsSS-LPR}.

The goal of specifying the sampling models is to indicate the exchangeability assumptions for the data.  To fully specify the necessary likelihood functions, $\theta_a$ will denote the parameter. For the common-source problem, the model selection is then a selection between $M_1$ and $M_2$ for the two sets of unknown-source observations.  Under $M_1$, the likelihood function for the unknown-source observations will be denoted 
\begin{equation}
\mathcal{L}(\theta_{a}|\mathbf{x_b}, \mathbf{x_c}, M_1) = f(\mathbf{x_b}, \mathbf{x_c}|\theta_{a}, M_1) = \int f(\mathbf{x_b}, \mathbf{x_c}|a_u,\theta_{a}) dG_{\theta_a}(a_u).
\end{equation} 
Under $M_2$, the likelihood function for the unknown-source observations will be denoted 
\begin{eqnarray}
\mathcal{L}(\theta_{a}|\mathbf{x_b}, \mathbf{x_c}, M_2) &=& \mathcal{L}(\theta_{a}|\mathbf{x_b}, M_2)\mathcal{L}(\theta_{a}| \mathbf{x_c}, M_2) = f(\mathbf{x_b}|\theta_{a}) f(\mathbf{x_c}|\theta_{a}) \notag
\\ &=& \int f(\mathbf{x_b}|a_b,\theta_{a}) dG_{\theta_a}(a_b) \int f(\mathbf{x_c}|a_c,\theta_{a}) dG_{\theta_a}(a_c).
\end{eqnarray} 
These likelihood functions will be used in Section~\ref{methods} to define statistics from both the classical and the Bayesian paradigms needed to perform the model selection.

\subsection{Specific Source Models}\label{ss}

\begin{figure} 
\includegraphics[width=\textwidth]{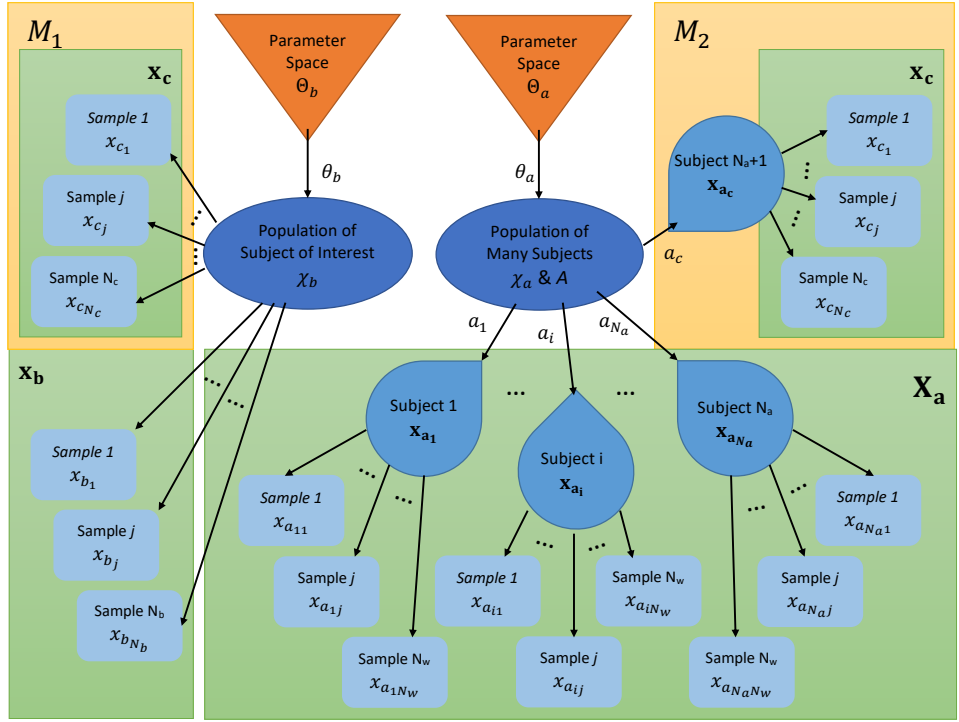}
\caption[]{Hierarchy and notation for the specific-source model selection problem where the green boxes denote the datasets and the yellow boxes denote the two competing models.}
\label{SSmodels}
\end{figure}

In contrast to the common-source problem, the idea of the specific-source non-nested model selection problem is that you have built up a dataset, $\mathbf{X_a}$, of observations from many different subjects in some population and there is one particular subject of interest (that is not a subject from the population mentioned previously) that may be the source of another set of observations. In addition to the dataset $\mathbf{X_a}$ (the number of subjects is denoted $N_a$, and the number of samples from within each subject is denoted $N_w$) you have also collected a set of observations from the particular subject of interest.  This dataset is denoted $\mathbf{x_b}$ (the sample size is $N_b$) and is composed of observations from a second, independent population, consisting of a single source.  Now, when a set of observations with an unknown source $\mathbf{x_c}$ (the sample size is $N_c$) is obtained, you want to determine whether this observation has come from the population associated with the subject of interest, denoted by $M_1$, or if it has been generated by a randomly selected subject from the larger population of many sources, denoted by $M_2$. The corresponding sampling models for the entire set of data, denoted $\mathbf{X_N} = \{ \mathbf{X_a},  \mathbf{x_b}, \mathbf{x_c}\}$ where $N$ denotes the index for the various sample sizes, are given in the supplementary material \citet{BFvsLRsupplement} and visualized in Figure~\ref{SSmodels}.  An example demonstrating the set-up of the specific-source problem for a forensic science application is given in Ommen and Saunders \citet{CSvsSS-LPR}. 

Similar to the common-source problem, the sampling models only provide information about the exchangeability of the observations so the parametric models need to be specified as well.  Let $\theta_a$ denote the parameter associated with the large population of many subjects and let $\theta_b$ denote the parameter associated with the specified subject's population. For the specific-source problem, the first model implies that $\mathbf{x_c}$ has been generated according to model for the specified subject, whereas the second model implies that $\mathbf{x_c}$ has been generated according to model for the large population of many subjects. Similar to the common-source problem, the model selection problem is then a selection between $M_1$ and $M_2$ for $\mathbf{x_c}$. Under the first model, the likelihood function for the unknown-source observation is denoted 
\begin{equation}
\mathcal{L}(\theta_b|\mathbf{x_c}, M_1) = f(\mathbf{x_c}|\theta_b, M_1) = \prod_{i=1}^{N_c}f(x_{c_i}|\theta_b)
\end{equation} 
and under the second model it will be denoted
\begin{equation}
\mathcal{L}(\theta_a|\mathbf{x_c}, M_2) = f(\mathbf{x_c}|\theta_a, M_2) = \int f(\mathbf{x_c}|a_c,\theta_a) dG_{\theta_a}(a_c).
\end{equation}  
These likelihood functions will be used in Section~\ref{methods} to define statistics from both the classical and the Bayesian paradigms needed to perform the model selection.

\subsection{Discussion}

While the non-nested model selection problems discussed in this article have been motivated by forensic science applications, they can be applied to a number of other areas as well.  The common-source and specific-source frameworks for non-nested model selection might be used in any areas where statistical pattern recognition methods can be applied.  For example, consider the situation in which there is a database of observations, $\mathbf{X}$, that have been divided into $k$ different classes, and you wish to classify a unknown-source set of observations, $\mathbf{y}$ to one of these $k$ classes.  Then the traditional Bayes' classifier is given by
\begin{equation*}
r(\mathbf{y}) = \argmax_{i \in \{1, 2, \ldots, k\}} P(M_i|\mathbf{y})
\end{equation*} 
where $P(M_i|\mathbf{y})$ is the frequentist posterior probability of the model for $i^{th}$ class creating the unknown-source set of observations $\mathbf{y}$  \citet{TWAnderson}.  Now, suppose that one of those $k$ classes is of particular interest to you for some reason.  The model for this class can be used in the specific-source framework to describe the generation of data from the population associated with the single subject of interest and the combination of all the other models for the other $k-1$ classes can be used in the specific source framework to describe the generation of data from the larger population of many subjects.  Then, the specific-source (non-nested) model selection problem would be a selection between the model that generated the unknown-source set of observations, $\mathbf{y}$.  Or, consider the scenario in which you have two unknown-source sets of observations, $\mathbf{y_1}$ and $\mathbf{y_2}$, and you want to know whether or not they come from the same class, without specifying which of the $k$ classes.  Then the combination of all the models for each of the $k$ classes serves as the model associated with the population for the common-source framework, and the (non-nested) model selection problem is then a decision about whether $\mathbf{y_1}$ and $\mathbf{y_2}$ have come from the same unspecified class or from two different unspecified classes.

The main difference between the model selection problems discussed in this paper and traditional classification or pattern recognition problems is the order in which you observe the sets of data.  In forensics applications, the unknown-source data is usually observed first, and then the other sets of data are subsequently collected.  This is because the unknown-source data is typically collected from the crime scene.  Investigators don't collect samples from a person of interest unless a crime has been committed first, so the samples from the known subjects are collected last.  In traditional classification problems, the known-source data is typically collected first, and then the unknown-source data is analyzed as it is observed. Note, the order in which the data is observed is not a necessary condition for application of this problem set-up. However, the results in this paper rely on the fixed observation of the unknown-source data, making them particularly useful for forensic science, and other similar, applications where this subset of the data is observed first.

\section{Methods}\label{methods}

In general, Bayesian statistical methods of approaching the forensic identification of source problems have been advocated by a number of forensic science researchers, especially in Europe \citet{ENFSI, AitkenLPR, SciJusBerger, SciJusBiedermann} while alternative methods have been advocated by a number of forensic science researchers primarily in the United States \citet{HariSteve, FRStat, Kafadar-noLR}.  The Bayesian methods are centered around computing the Bayes Factor which provides a relative measure for how much the data supports the two competing models.  In contrast, the alternative methods typically focus on finding an approximation to the likelihood ratio.  One of the most confusing things about these problems in forensics is that both of these values are referred to as the ``likelihood ratio" regardless of which method was used to compute the value.  It is very easy to imagine that a similar phenomena occurs in related areas outside of forensics, too.  Therefore, it is very important to distinguish the methods used to compute these values. The following subsections will provide a very brief review of each of the methods and the notation to define the likelihood ratio and the Bayes Factor specifically for each of the non-nested model selection problems we are interested in.

\subsection{The Likelihood Ratio}\label{LR-section}

Typically, the term likelihood ratio is used in reference to the likelihood ratio test statistic developed by Neyman and Pearson for simple- and nested-model selection.  Another common reference for the term likelihood ratio is the log-likelihood ratio statistic for nested-model selection. However, neither of these methods are directly applicable to non-nested models, like the two model selection problems introduced above.  Since these early developments, others have explored the use of the likelihood ratio in more generality for non-nested model selection (see Vuong \citet{Vuong} for example).  In this article, the phrase likelihood ratio function will be used to describe the ratio of the two different likelihood functions for the unknown-source observations under the two competing models.  In particular, the likelihood ratio function for the common-source model selection problem described in Section~\ref{cs} is given by
\begin{eqnarray}\label{LRcs}
\lambda_{cs}(\theta_a| \mathbf{x_b}, \mathbf{x_c}) &=& \dfrac{\mathcal{L}(\theta_{a}|\mathbf{x_b}, \mathbf{x_c}, M_1)}{\mathcal{L}(\theta_{a}|\mathbf{x_b}, M_2)\mathcal{L}(\theta_{a}| \mathbf{x_c}, M_2)}
\end{eqnarray}
and the likelihood ratio function for the specific-source model selection problem described in Section~\ref{ss} is given by
\begin{equation}\label{LRss}
	\lambda_{ss}(\theta_{a},\theta_{b}| \mathbf{x_c}) = \dfrac{\mathcal{L}(\theta_{b}|\mathbf{x_c}, M_1)}{\mathcal{L}(\theta_{a}|\mathbf{x_c}, M_2)} .
\end{equation}
For given observations of the unknown-source data, the likelihood ratio function for the common-source model selection problem is a function of the multivariate parameter vector $\theta_a$ which is constrained to take values in the space $\Theta_a$.  In the specific-source model selection problem, the likelihood ratio function is a multivariate function of the joint parameter vector for both $\theta_a$ and $\theta_b$, given the observation of the unknown-source data. In this sense, the value of the likelihood ratio function which corresponds to the values of the parameters associated with the true sampling distributions implied by models $M_1$ and $M_2$ can be considered the ``true likelihood ratio."  The parameter value $\theta_{a_0}$ is used to indicate the true parameter for the larger population of many sources while $\theta_{b_0}$ is used to indicate the true parameter for the population associated with the single source of interest. Since we do not know what the values of these parameters are, the true likelihood ratio is fixed, but unknown.  The form of the true likelihood ratio for the common-source model selection problem is given by
\begin{eqnarray}\label{trueLRcs}
	\lambda_{cs}(\theta_{a_0}| \mathbf{x_b}, \mathbf{x_c}) &=& \dfrac{\mathcal{L}(\theta_{a_0}|\mathbf{x_b}, \mathbf{x_c}, M_1)}{\mathcal{L}(\theta_{a_0}|\mathbf{x_b}, M_2)\mathcal{L}(\theta_{a_0}| \mathbf{x_c}, M_2)}  
\end{eqnarray}
and for the specific-source model selection problem by
\begin{equation}\label{trueLRss}
	\lambda_{ss}(\theta_{a_0},\theta_{b_0}| \mathbf{x_c}) = \dfrac{\mathcal{L}(\theta_{b_0}|\mathbf{x_c}, M_1)}{\mathcal{L}(\theta_{a_0}|\mathbf{x_c}, M_2)}.
\end{equation}
Note that the true likelihood ratio given by either Equation~\ref{trueLRcs} or Equation~\ref{trueLRss} is not the limiting form of the likelihood ratio function as the size of entire set of data (including the size of the observations with unknown source) grows infinitely larger.  Theorem 5.1 from Vuong \citet{Vuong} provides the limiting form of the likelihood ratio function for this scenario, which will not be considered further in this article.  The true likelihood ratio will take values between zero and positive infinity (exclusive of the endpoints) since the size of the unknown-source observations remains fixed.

\subsection{The Bayes Factor}\label{BF-section}

The Bayes Factor is a well-studied statistic often used for performing model selection on a wide variety of models \citet{KassRaf, NewRaf}. One of the advantages of the using the Bayes Factor for model selection is that the models do not need to be nested.  The Bayes Factor can be viewed as the strength of the support the data provides for the two competing models.  A Bayes Factor greater than one indicates that the data support model $M_1$ over model $M_2$, in contrast to a Bayes Factor less than one which indicates that the data support model $M_2$ over model $M_1$.  A Bayes Factor that is equal to one means that the data cannot discriminate between the two competing models.  However, the cost for increased applicability of the Bayes Factor for model selection is the need to specify prior probabilities for each of the models.  The Bayes Factor can then be used to update the prior odds for the models to arrive at the posterior odds.  

\begin{equation}\label{BF}
\underbrace{
\frac
{\mathbb{P}\left( M_1 | \mathbf{X} \right)}
{\mathbb{P}\left( M_2 | \mathbf{X} \right)}
}_\text{Posterior Odds}  
=
\underbrace{
\frac
{\mathbb{P}\left( \mathbf{X} | M_1 \right)}
{\mathbb{P}\left( \mathbf{X} | M_2 \right)}
}_{
\begin{smallmatrix}
   \text{Bayes Factor}
\end{smallmatrix}
}
\times \ 
\underbrace{
\frac
{\mathbb{P}\left( M_1 \right)}
{\mathbb{P}\left( M_2 \right)}
}_\text{Prior Odds} 
\end{equation}
The posterior odds will then be used to select the model.  A posterior odds greater than one indicates that, given all the data you have observed, model $M_1$ is preferred to model $M_2$, whereas a prior odds less than one indicates that model $M_2$ is preferred to model $M_1$.  

In the situation where the models involve parametric distributions, the Bayes Factor is generally given by the following expression:
\begin{equation}\label{genV}
\beta(\mathbf{X}) = \dfrac{\int f(\mathbf{X}| \theta, M_1)\ d\Pi(\theta|M_1)}{\int f(\mathbf{X}| \theta, M_2)\ d\Pi(\theta|M_2)}.
\end{equation}
This can be interpreted as a ratio of marginal (or sometimes called integrated) likelihoods for the two competing models.  Notice that the Bayes Factor itself contains a second set of prior distributions, denoted by $\Pi$, for the values of the parameters under each of the two competing models.  These prior distributions are distinct from the prior odds for the models.

The precise form of the Bayes Factor for the common-source model selection problem described in Section~\ref{cs} is derived under the assumption that the prior distribution of $\theta_a$ will be the same under both models, and is given by
\begin{eqnarray}\label{Vcs_fact1}
\beta_{cs}(\mathbf{X_N}) &=& \dfrac{\int f(\mathbf{x_b}, \mathbf{x_c}| \theta_a, M_1)\ d\Pi(\theta_a|\mathbf{X_a})}{\int f(\mathbf{x_b}| \theta_a, M_2) f(\mathbf{x_c}| \theta_a, M_2)\ d\Pi(\theta_a|\mathbf{X_a})}
\\ &\equiv & \dfrac{m(\mathbf{x_b}, \mathbf{x_c}|\mathbf{X_a}, M_1)}{m(\mathbf{x_b}, \mathbf{x_c}|\mathbf{X_a}, M_2)}. \notag
\end{eqnarray}
Under the assumption that the prior distribution of $\theta_a$ is statistically independent of the prior distribution for $\theta_b$, the Bayes Factor for the specific-source model selection problem described in Section~\ref{ss} is given by
\begin{eqnarray}\label{Vss_fact1}
\beta_{ss}(\mathbf{X_N}) &=& \dfrac{\int f(\mathbf{x_c}|\theta_{b}, M_1)\ d\Pi(\theta_b|\mathbf{x_b})}{\int f(\mathbf{x_c}|\theta_{a}, M_2)\ d\Pi(\theta_a|\mathbf{X_a})}
\\ &\equiv & \dfrac{m(\mathbf{x_c}|\mathbf{x_b}, M_1)}{m(\mathbf{x_c}|\mathbf{X_a}, M_2)}. \notag
\end{eqnarray}
The derivations of Equation~\ref{Vcs_fact1} and Equation~\ref{Vss_fact1} are given in the appendix of Ommen et al. 2017 \citet{dmo_mcse}.
Bayes Factors in these forms are often very computationally intensive to compute since the integrals rarely have closed-form solutions \citet{KassRaf, NewRaf, dmo_mcse}. 

\section{Relationships between the LR and the BF}\label{LRvsBF1}

As we previously mentioned, the likelihood ratio is most often used within the classical paradigm of statistics, while the Bayes Factor is most often associated with the subjective Bayesian paradigm.  In this section, we will provide some expressions which directly relate these statistics from the two different paradigms.  These expressions will serve as the foundation for subscribers of one statistical paradigm to more effectively communicate their results to subscribers of differing paradigms.  For example, imagine that a statistician working within the classical paradigm presents a model selection result, in the form of the likelihood ratio function, to a statistician working within the Bayesian paradigm.  Then, the corresponding Bayes Factor for the model selection problem can be expressed as the expected value of the likelihood ratio function with respect to the Bayesian's subjective posterior belief about the parameters of the sampling models given the entire collection of data under the second model.  This is formalized in the following two equations below for the common-source and specific-source model selection problems, respectively.  The derivations of Equation~\ref{Vcs_fact2} and Equation~\ref{Vss_fact2} can be found in Appendix~\ref{A-derive}.

\begin{eqnarray}\label{Vcs_fact2}
\beta_{cs}(\mathbf{X_N}) &=& \int \lambda_{cs}(\theta_a| \mathbf{x_b}, \mathbf{x_c})\ d\Pi(\theta_a|\mathbf{X_N}, M_2)
\end{eqnarray}

\begin{eqnarray}\label{Vss_fact2}
\beta_{ss}(\mathbf{X_N}) &=& \int \lambda_{ss}(\theta_{a},\theta_{b}| \mathbf{x_c})\ d\Pi(\theta_{a},\theta_{b}|\mathbf{X_N}, M_2)
\end{eqnarray}

A similar expression can be derived using the first model for the data in the posterior distribution, with the introduction of two inverses.  The derivations of Equation~\ref{Vcs_fact2p} and Equation~\ref{Vss_fact2p} can also be found in Appendix~\ref{A-derive}.

\begin{equation}\label{Vcs_fact2p}
\beta_{cs}(\mathbf{X_N}) = \left[ \int \dfrac{1}{\lambda_{cs}(\theta_a| \mathbf{x_b}, \mathbf{x_c})}\ d\Pi(\theta_a|\mathbf{X_N}, M_1) \right]^{-1}
\end{equation}

\begin{equation}\label{Vss_fact2p}
\beta_{ss}(\mathbf{X_N}) = \left[ \int \dfrac{1}{\lambda_{ss}(\theta_{a},\theta_{b}| \mathbf{x_c})}\ d\Pi(\theta_{a},\theta_{b}|\mathbf{X_N}, M_1) \right]^{-1}
\end{equation}

In these expressions, notice that the data with unknown source is included in the posterior distribution for the parameters.  Another interesting thing to note about these forms for the Bayes Factor is that one of them will be computed using a misspecified model for the data $\mathbf{X_N}$ depending on the truth of which model actually generated the evidence.  

\subsection{Asymptotic Results}\label{LRvsBF2}
In the reviewed literature, asymptotic properties of the Bayes Factor, particularly in the case of non-nested model selection, have been examined with respect to data from a single population representing a single source of information when the number of observations from that single population is growing.  In this case, the Bayes Factor will diverge to positive infinity (in probability) when the model in the numerator is preferred and will converge to zero (in probability) when the model in the denominator is preferred \citet{Chib-BFconsistency}. In the case of both the common-source and specific-source model selection problems, there are multiple sources of information and only the number of observations from a portion of these sources is allowed to grow.  In this section, we examine the consistency of the Bayes Factor by way of a well-known result, the Bernstein-von Mises Theorem \citet{VDV}.  The Bernstein-von Mises Theorem is reproduced from Van der Vaart \citet{VDV} for clarity in the supplementary material \citet{BFvsLRsupplement}.  We will show that, under certain regularity conditions, the Bayes Factor for both the common-source and specific-source model selection problems will converge to the true likelihood ratio. These results are formalized in the theorems to follow. The notational modifications used to facilitate the asymptotic results and corresponding proofs are given in Appendix~\ref{B-proof}.

\begin{quote}\normalsize
\begin{thm}[Common-Source Bayes Factor Consistency]\label{CS-BF-LR}
Given a fixed observation of $\mathbf{x_b}$ and $\mathbf{x_c}$, suppose that the likelihood ratio function $\lambda_{cs}(\theta_a|\mathbf{x_b}, \mathbf{x_c})$ is bounded in a neighborhood of $\theta_{a_0}$ and that $\hat{\theta}_{a|M_2}$ is a consistent estimator for $\theta_{a_0}$.  Furthermore, suppose the assumptions of the Bernstein-von Mises Theorem are satisfied.  Then as $n_a \to \infty$, the Bayes Factor converges in $P_{\theta_a}^{n_a}$-probability to the true likelihood ratio,
\begin{equation*}
\beta_{cs}({X_{a,n_a}}, \mathbf{x_b}, \mathbf{x_c})  \overset{P_{\theta_a}^{n_a}}{\longrightarrow}  \lambda_{cs}(\theta_{a_0}|\mathbf{x_b}, \mathbf{x_c}).
\end{equation*}
\end{thm}
\end{quote}
The proof of this theorem is provided in Appendix~\ref{B-proof}. 

\begin{quote}\normalsize
\begin{thm}[Specific-Source Bayes Factor Consistency]\label{SS-BF-LR}
Given a fixed observation of $\mathbf{x_c}$,  suppose that the likelihood ratio function, $\lambda_{ss}(\theta| \mathbf{x_c})$, is bounded in a neighborhood of $\theta_0$ and that $\hat{\theta}_{n|M_2}$ is a consistent estimator for $\theta_0$.  Furthermore, suppose the assumptions of the Bernstein-von Mises Theorem are satisfied. Then the Bayes Factor converges in $P_{\theta}^n$-probability to the true likelihood ratio as $n \to \infty$,
\begin{equation*}
\beta_{ss}({X_{a,n_a}}, X_{b,n_b}, \mathbf{x_c})  \overset{P_{\theta}^n}{\longrightarrow}  \lambda_{ss}(\theta_0| \mathbf{x_c}).
\end{equation*}
\end{thm}
\end{quote}
The proof of this theorem is provided in Appendix~\ref{B-proof}. 

These results expand upon the foundation for members of one statistical paradigm to communicate more effectively with members of another paradigm.  For example, imagine that a statistician working within the Bayesian paradigm presents a model selection result, in the form of the Bayes Factor, to a statistician working within the classical paradigm.  Then, the likelihood ratio for the corresponding model selection problem can be thought of as the limit of the given Bayes Factor if more and more information is gathered from a subset of the data sources.  Unfortunately, this is not a very useful result in practice since it is often infeasible to gather more information, and we rarely have large enough samples sizes for the value of the Bayes Factor to be used as a large-sample replacement for the true value of the likelihood ratio.

\subsection{Bayesian Credible Intervals for the LR}\label{LRvsBF3}
Since the results from the previous section only provide theoretical information for transforming the Bayes Factor into the corresponding likelihood ratio, the purpose of this section is to provide statisticians with a practical way to relate the two.  Credible intervals, particularly used to determine posterior concentration rates, are a popularly researched topic recently in the high-dimensional and non-parametric statistics literature (see Hoffmann et al. \citet{Rousseau2015}, Rockova and van der Pas \citet{Rockova-trees}, Donnet et al. \citet{Rousseau2018}, and Nickl and Sohl \citet{Sohl} for example).  Using similar methods, we will create a credible interval for the value of the likelihood ratio derived from the posterior distribution for the parameters.  The necessary notational modifications related to the result and corresponding proofs are given in Appendix~\ref{C-ints}.

\begin{quote}\normalsize
\begin{thm}[Approximate $1-\alpha$ Credible Interval for the LR]\label{CIexist}
Let let the assumptions of Lemma~\ref{C-ints-1} and Lemma~\ref{C-ints-2} hold and let 
$$\mathcal{I}_n = \beta({X_n}) \pm \Phi^{-1}(1-\alpha/2)\sigma_n$$
where $\beta({X_n})$ represents the sequence of either the common-source or the specific-source Bayes Factor, $0<\alpha<1$ is the desired significance level, $\sigma_n$ is the sequence of posterior standard errors for the likelihood ratio, and $\Phi^{-1}$ is the standard normal quantile function. Then as $n \to \infty$
\begin{equation*}
	\Pi(\lambda(\theta) \in \mathcal{I}_n |{X_n}, M_2) \to 1-\alpha.
\end{equation*}
\end{thm}	
\end{quote}
Please see Appendix~\ref{C-ints} for Lemma~\ref{C-ints-1} and Lemma~\ref{C-ints-2} and for a proof of this theorem. We would like to note that we chose a form of the interval that is guaranteed to contain the posterior mean, in this case the Bayes Factor, instead of choosing an equal-tails or highest posterior density. In our experience, the posterior distributions may be so skewed that the Bayes Factor doesn't actually fall in the body of the posterior distribution, an therefore may not be contained in the interval constructed using these other methods. This would have serious implications in applications like forensic science.

The credible interval described above provides a range of probable likelihood ratio values that would correspond to the subjective, personal Bayes Factor provided for the non-nested model selection problem.  Now, this method has a couple of drawbacks.  The first is that this interval for a classical statistic relies on subjective Bayesian probabilities.  The interval is highly dependent on the given value of the (possibly estimated) Bayes Factor and the corresponding prior distributions chosen.  Another disadvantage of this method is that if the prior has not been chosen properly, the credible interval for the LR corresponding to the given Bayes Factor may \textit{not} actually contain the true value of the likelihood ratio.  In this way, the credible interval is easily misinterpreted.  The interval must not be interpreted as a range of probable values for the true likelihood ratio (as if you had full, infinite data for every source of information), but instead should be interpreted as a range of probable values for the estimated likelihood ratio given the limited data and given the chosen prior distribution.  This interpretation may be unsatisfying for the classical statistician who must rely on someone else's subjective belief to determine a corresponding value of the estimated likelihood ratio.  Finally, if the value of the Bayes Factor has been estimated, using Monte Carlo integration for example, it is unclear how, or even if, the computational error should be incorporated into this credible interval for the corresponding LR.

\section{Example}

\begin{figure} 
\includegraphics[width=\textwidth]{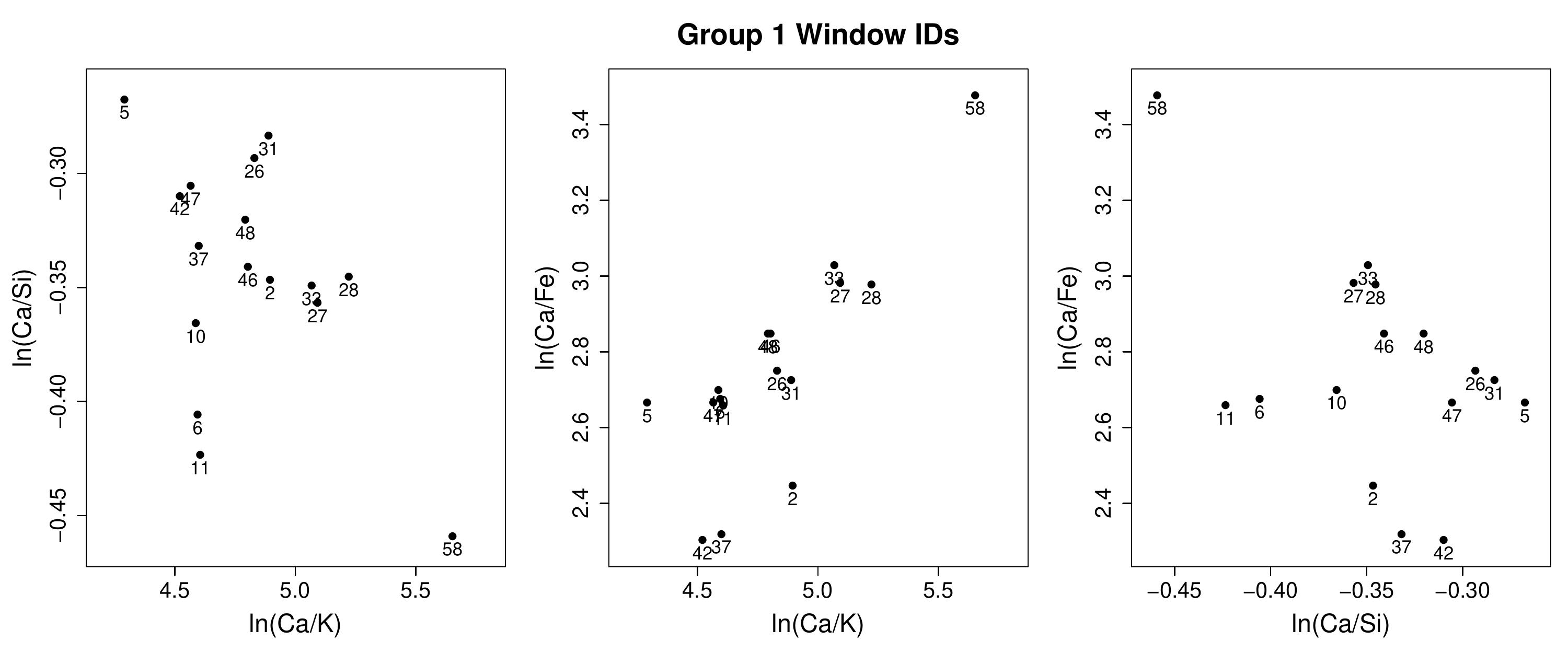}
\caption[]{Pairwise plots of the mean elemental concentrations for each window in the first group along with the corresponding window identification number.}
\label{glassIDs}
\end{figure}

\begin{figure} 
\includegraphics[width=\textwidth]{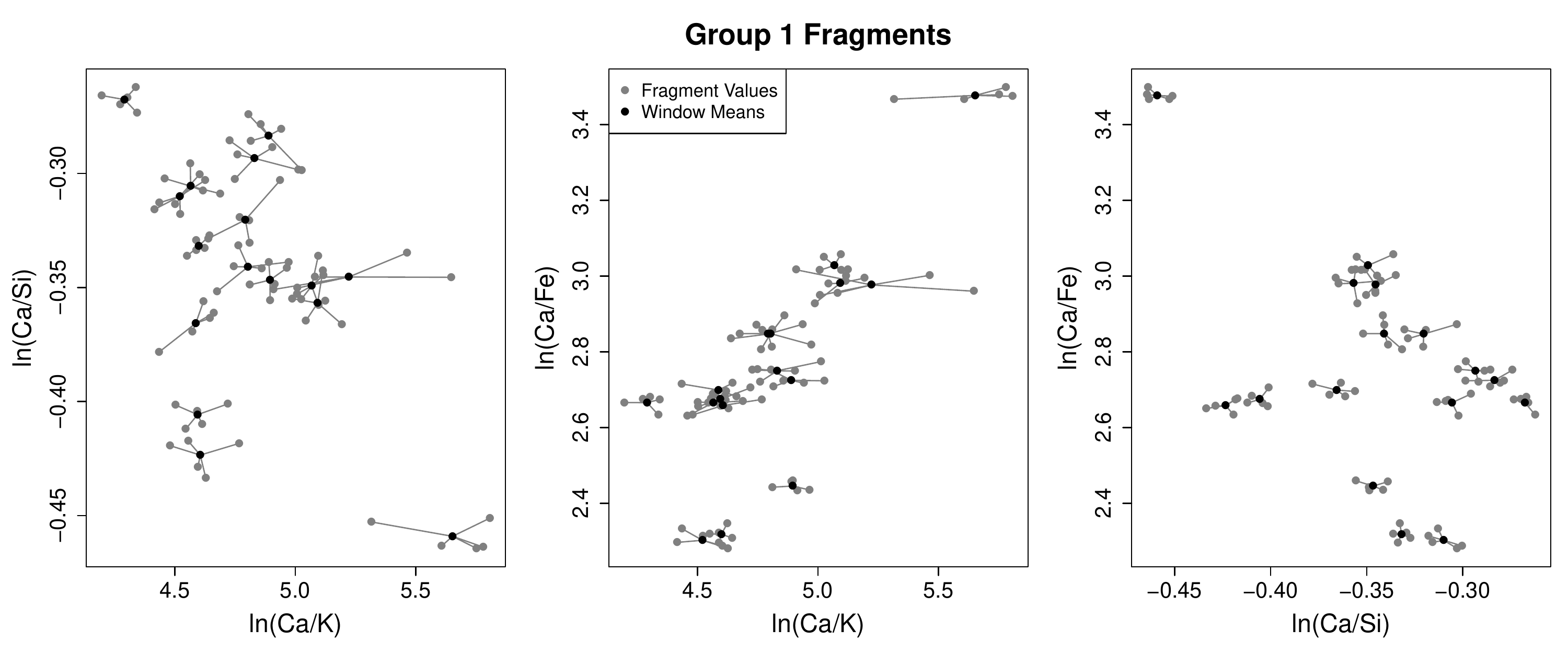}
\caption[]{Pairwise plots of the elemental concentrations for each fragment in the first group along with the mean elemental concentrations for each associated window.}
\label{glass}
\end{figure}

One of the areas of forensic science where both the common-source and specific-source frameworks would straightforwardly apply is trace elemental analysis of glass evidence.  Ommen et al. \citet{dmo_mcse} provides a set-up for both the common-source and specific-source Bayes Factors for the trace elemental composition of a dataset of 62 different window panes, each with 5 measured fragments of glass.  This data set was originally collected by Dr. JoAnn Buscaglia of the Federal Bureau of Investigation Laboratory Division and analyzed by Aitken and Lucy using a multivariate analysis \citet{AitkenLucy}.  For this example, the 16 window panes from the first group will be used as the data, $\mathbf{X_N}$.  Figure~\ref{glassIDs} provides the window ID numbers corresponding to the pairwise values for the transformed trace elemental compositions provided in Figure~\ref{glass}.

Directly computing the value of the likelihood ratio for this dataset for either the common-source or specific source frameworks, given by Equation~\ref{LRcs} and Equation~\ref{LRss}, respectively, is impossible since we do not know the parameter values under either model.  Also, approximating the value of the true likelihood ratio using the asymptotic results given by Theorem~\ref{CS-BF-LR} or Theorem~\ref{SS-BF-LR} is inadvisable since we have such small sample sizes.  So, this would be a case in which we would need to compute credible intervals for the likelihood ratio using Theorem~\ref{CIexist} in order to compare the results from the two differing statistical paradigms.  The Bayes Factors needed to compute the intervals for the common-source problem are given by Equation~\ref{Vcs_fact2} and Equation~\ref{Vcs_fact2p} or by Equation~\ref{Vss_fact2} and Equation~\ref{Vss_fact2p} for the specific-source problem.  Ommen et al. gives the Bayes Factors in the form of Equation~\ref{Vcs_fact1}, Equation~\ref{Vss_fact1}, and Equation~\ref{Vcs_fact2}, but the remaining three required forms follow readily.  

In the interest of space, we will only consider the common-source example for this small glass dataset.  Following the setup given in Ommen et al. \citet{dmo_mcse}, the transformed trace elemental concentrations will be modeled by the same normal distributions and using the same normal prior distributions for the mean parameters and the same inverse Wishart prior distributions for the covariance parameters. Table~\ref{example1} gives the value of the Bayes Factor, posterior standard deviation, and the corresponding credible interval for the likelihood ratio for a sampling of window panes serving as $\mathbf{x_b}$, with window pane 10 serving as $\mathbf{x_c}$, and using the second group of window panes to determine the prior hyperparameters. Table~\ref{example2} gives the value of the Bayes Factor, posterior standard deviation, and the corresponding credible interval for the likelihood ratio for the same sampling of window pairs serving as the data for the unknown-source observations, $\mathbf{x_b}$ and $\mathbf{x_c}$, but using the third group of window panes to determine the prior hyperparameters.  This would be like having two different experts compute the same values using two different prior distributions. In both examples, the remaining 14 window panes are used for $\mathbf{X_a}$. Details of the example are given in the supplementary material \citet{BFvsLRsupplement}.

\begin{table}
\centering
\caption{Bayes Factors, posterior standard deviations, and corresponding 95\% credible intervals\\ truncated at zero for the likelihood ratio from Expert 1}\label{example1}
\begin{tabular}{l|cccccccccccccccc}
\hline
&\multicolumn{3}{c}{$\mathbf{x_b}$} \\
& 10 & 6 & 48 \\
\hline
Bayes Factor & 779.30 & $5.10\times10^{-6}$ & $3.04\times10^{-10}$ \\
Posterior Std. Dev.& 249.7349 & $9.84\times10^{-5}$ & $4.35\times10^{-9}$ \\
Credible Interval & (289.83, 1268.77) & (0, $1.98\times10^{-4}$) & (0, $8.82\times10^{-9}$) \\
\hline
\end{tabular}
\end{table}

\begin{table}
\centering
\caption{Bayes Factors, posterior standard deviations, and corresponding 95\% credible intervals\\ truncated at zero for the likelihood ratio from Expert 2}\label{example2}
\begin{tabular}{l|cccccccccccccccc}
\hline
&\multicolumn{3}{c}{$\mathbf{x_b}$} \\
& 10 & 6 & 48 \\
\hline
Bayes Factor & 186.46 & $7.11\times10^{-7}$ & $4.54\times10^{-9}$ \\
Posterior Std. Dev. & 48.8764 & $8.86\times10^{-6}$ & $8.41\times10^{-8}$ \\
Credible Interval & (87.66, 279.25) & (0, $1.81\times10^{-5}$) & (0, $1.69\times10^{-7}$) \\
\hline
\end{tabular}
\end{table}

As you can see from Table~\ref{example1} and Table~\ref{example2}, the credible intervals for the likelihood ratios have been truncated at zero.  We did this because it is a well-known fact that likelihood ratios cannot be negative.  If we did not truncate the lower endpoints at 0, then the endpoints would be negative.  This is due to the very small sample sizes that we are working with.  In addition, it is possible using this method to have the credible intervals for the likelihood ratio from two different experts be non-overlapping.  This means that Expert 1 and Expert 2 would likely disagree on the value of evidence given the limited size of the datasets due to the fact that they used different prior distributions. However, it is worth noting that neither of the intervals contain the neutral value of one, and so neither would present misleading evidence to a fact-finder.  In the event that more data was collected, the intervals from the two experts would eventually overcome the disagreement caused by the differing priors. This is a direct result from Theorem~\ref{CIexist}.  Exactly how much additional data would need to be collected to overcome the disagreement in a practical scenario may be estimated through simulation (although this is not directly considered in this example).

In this example, the glass dataset was collected for research purposes, and not intended to represent the type of evidence forensic scientists would see in casework.  However, the two non-nested models provided in Ommen et al. \citet{dmo_mcse} can be applied to trace elemental analysis of glass evidence in general.  In addition, the non-nested models can be modified to work for trace elemental analysis of any type of evidence, not just glass.  For instance, the common-source and specific-source frameworks can be used for trace elemental analysis of ink \citet{ICPMS-ink}, paper\citet{ICPMS-paper}, soil \citet{ICPMS-soil}, and paint \citet{ICPMS-paint} as well.  Even more generally, these two frameworks can be applied to evidence measured by nearly any quantitative method.  Once the appropriate framework has been set, then all of the results relating the Bayes Factor to the likelihood ratio discussed in this article apply.

\section{Conclusions}
In this article, we have expressed the Bayes Factor for two non-nested model selection problems as the expected value of the corresponding likelihood ratio function with respect to the posterior distribution for the parameters given the entire set of data where the set(s) of unknown-source observations has been generated according to the second model. This expression has led to a number of useful theoretic and practical results relating the Bayesian solution of these non-nested model selection problems to the classical solution.  These non-nested model selection problems differ from those previously considered in other asymptotic results regarding the Bayes Factor in that more than one source of information is generating data.  In addition, for the asymptotic results, only a portion of these datasets are allowed to grow in size while some observations remained fixed throughout.  These results have strong implications for the field of forensic science where the term ``likelihood ratio" is used ubiquitously to mean a solution to the forensic identification of source problem, regardless of the style of statistics used in the analysis.

\section*{Acknowledgements}
We would like to thank Dr. JoAnn Buscaglia for providing the dataset used in the example and for the many helpful discussions throughout the years. Since the research presented in this article was performed as part of Dr. Ommen's dissertation research at South Dakota State University, she would like to thank the members of her departmental Ph.D. committee, Dr. Kurt Cogswell and Dr. Cedric Neumann, for their support and guidance.

\bibliography{DanicaPhD_Bib}
\bibliographystyle{plainnat}

\appendix

\section{Bayes Factor Derivations}\label{A-derive}

Further details of all the derivations that follow are given in the supplementary material \citet{BFvsLRsupplement}.

\subsection{Equation~\ref{Vcs_fact2}}
This derivation is a generalization of the derivation given in Ommen et al.\citet{dmo_mcse}.

\begin{subequations}
\begin{eqnarray}
\beta_{cs}(\mathbf{X_N}) \hspace{-5pt}&=&  \mathlarger{\int} \dfrac{f(\mathbf{x_b}, \mathbf{x_c}| \theta_a, M_1)}{m(\mathbf{x_b}, \mathbf{x_c}|\mathbf{X_a}, M_2)}\ d\Pi(\theta_a|\mathbf{X_a})
\\ &=& \mathlarger{\int} \dfrac{f(\mathbf{x_b}, \mathbf{x_c}| \theta_a, M_1)}{f(\mathbf{x_b}, \mathbf{x_c}| \theta_a, M_2)} \times \dfrac{f(\mathbf{x_b}, \mathbf{x_c}| \theta_a, M_2)}{m(\mathbf{x_b}, \mathbf{x_c}|\mathbf{X_a}, M_2)}\ d\Pi(\theta_a|\mathbf{X_a})
\\ &=& \mathlarger{\int} \dfrac{f(\mathbf{x_b}, \mathbf{x_c}| \theta_a, M_1)}{f(\mathbf{x_b}, \mathbf{x_c}| \theta_a, M_2)}\ d\Pi(\theta_a|\mathbf{X_a}, \mathbf{x_b}, \mathbf{x_c}, M_2)
\\ &=& \mathlarger{\int} \lambda_{cs}(\theta_a| \mathbf{x_b}, \mathbf{x_c})\ d\Pi(\theta_a|\mathbf{X_N}, M_2)
\end{eqnarray}
\end{subequations}

\subsection{Equation~\ref{Vss_fact2}}

This derivation closely follows the derivation of Equation~\ref{Vcs_fact2} given above.  
\begin{subequations}
\begin{eqnarray}
&\beta_{ss}(\mathbf{X_N})& \notag
\\ &=& \int \dfrac{f(\mathbf{X_a}| \theta_a) f(\mathbf{x_b}| \theta_b) f(\mathbf{x_c}|\theta_b, M_1)}{m(\mathbf{X_a}, \mathbf{x_b}, \mathbf{x_c}|M_2)}\ d\Pi(\theta_a,\theta_b)
\\ &=& \int \dfrac{f(\mathbf{x_c}|\theta_b, M_1)}{f(\mathbf{x_c}|\theta_a, M_2)} \times \dfrac{f(\mathbf{X_a}| \theta_a) f(\mathbf{x_b}| \theta_b) f(\mathbf{x_c}|\theta_a, M_2)}{m(\mathbf{X_a}, \mathbf{x_b}, \mathbf{x_c}|M_2)}\ d\Pi(\theta_a,\theta_b)
\\ &=& \int \dfrac{f(\mathbf{x_c}|\theta_b, M_1)}{f(\mathbf{x_c}|\theta_a, M_2)} \ d\Pi(\theta_a,\theta_b|\mathbf{X_a}, \mathbf{x_b}, \mathbf{x_c}, M_2)
\\ &=& \int \lambda_{ss}(\theta_a, \theta_b| \mathbf{x_c})\ d\Pi(\theta_a,\theta_b|\mathbf{X_N}, M_2)
\end{eqnarray}
\end{subequations}

\subsection{Equation~\ref{Vcs_fact2p}}

First, consider the reciprocal of the common-source Bayes Factor. 
\begin{subequations}
\begin{eqnarray}
\dfrac{1}{\beta_{cs}(\mathbf{X_N})} &=&  \mathlarger{\int} \dfrac{f(\mathbf{x_b}, \mathbf{x_c}| \theta_a, M_2)}{m(\mathbf{x_b}, \mathbf{x_c}|\mathbf{X_a}, M_1)}\ d\Pi(\theta_a|\mathbf{X_a})
\\ &=& \mathlarger{\int} \dfrac{f(\mathbf{x_b}, \mathbf{x_c}| \theta_a, M_2)}{f(\mathbf{x_b}, \mathbf{x_c}| \theta_a, M_1)} \times \dfrac{f(\mathbf{x_b}, \mathbf{x_c}| \theta_a, M_1)}{m(\mathbf{x_b}, \mathbf{x_c}|\mathbf{X_a}, M_1)}\ d\Pi(\theta_a|\mathbf{X_a})
\\ &=& \mathlarger{\int} \dfrac{f(\mathbf{x_b}, \mathbf{x_c}| \theta_a, M_2)}{f(\mathbf{x_b}, \mathbf{x_c}| \theta_a, M_1)}\ d\Pi(\theta_a|\mathbf{X_a}, \mathbf{x_b}, \mathbf{x_c}, M_1)
\\ &=& \mathlarger{\int} \dfrac{1}{\lambda_{cs}(\theta_a| \mathbf{x_b}, \mathbf{x_c})}\ d\Pi(\theta_a|\mathbf{X_N}, M_1)
\end{eqnarray}	
\end{subequations}
Therefore,
\begin{equation}
\beta_{cs}(\mathbf{X_N}) = \left[ \int \dfrac{1}{\lambda_{cs}(\theta_a| \mathbf{x_b}, \mathbf{x_c})}\ d\Pi(\theta_a|\mathbf{X_N}, M_1) \right]^{-1}	
\end{equation}

\subsection{Equation~\ref{Vss_fact2p}}

Similar to the common-source derivation above, consider the reciprocal of the specific-source Bayes Factor. 
\begin{subequations}
\begin{eqnarray}
&\dfrac{1}{\beta_{ss}(\mathbf{X_N})} & \notag
\\ &=& \int \dfrac{f(\mathbf{X_a}| \theta_a) f(\mathbf{x_b}| \theta_b) f(\mathbf{x_c}|\theta_a, M_2)}{m(\mathbf{X_a}, \mathbf{x_b}, \mathbf{x_c}|M_1)}\ d\Pi(\theta_a,\theta_b)
\\ &=& \int \dfrac{f(\mathbf{x_c}|\theta_a, M_2)}{f(\mathbf{x_c}|\theta_b, M_1)} \times \dfrac{f(\mathbf{X_a}| \theta_a) f(\mathbf{x_b}| \theta_b) f(\mathbf{x_c}|\theta_b, M_1)}{m(\mathbf{X_a}, \mathbf{x_b}, \mathbf{x_c}|M_1)}\ d\Pi(\theta_a,\theta_b)
\\ &=& \int \dfrac{f(\mathbf{x_c}|\theta_a, M_2)}{f(\mathbf{x_c}|\theta_b, M_1)} \times \dfrac{f(\mathbf{X_a}, \mathbf{x_b}, \mathbf{x_c}|\theta_a, \theta_b, M_1)}{m(\mathbf{X_a}, \mathbf{x_b}, \mathbf{x_c}|M_1)}\ d\Pi(\theta_a,\theta_b)
\\ &=& \int \dfrac{f(\mathbf{x_c}|\theta_a, M_2)}{f(\mathbf{x_c}|\theta_b, M_1)} \ d\Pi(\theta_a,\theta_b|\mathbf{X_a}, \mathbf{x_b}, \mathbf{x_c}, M_1)
\\ &=& \int \dfrac{1}{\lambda_{ss}(\theta_{a},\theta_{b}| \mathbf{x_c})}\ d\Pi(\theta_a,\theta_b|\mathbf{X_N},  M_1)
\end{eqnarray}
\end{subequations}
Therefore,
\begin{equation}
\beta_{ss}(\mathbf{X_N}) = \left[ \mathlarger{\int} \dfrac{1}{\lambda_{ss}(\theta_{a},\theta_{b}| \mathbf{x_c})}\ d\Pi(\theta_a,\theta_b|\mathbf{X_N}, M_1) \right]^{-1}	
\end{equation}

\section{Consistency Proofs}\label{B-proof}

The Bernstein-von Mises Theorem is a result which describes the contraction of the posterior distribution for the parameters as you observe more and more relevant data. The Bernstein-von Mises Theorem will be used to show that the Bayes Factor for the two different non-nested model selection problems is consistent towards the corresponding true likelihood ratio.  A version of the Bernstein-von Mises Theorem from van der Vaart \citet{VDV} is reproduced in the supplementary material \citet{BFvsLRsupplement} for ease of reference.  It should be noted that further details of the proofs in this section are also given in the supplementary material \citet{BFvsLRsupplement}.

Please observe the following notational extensions to facilitate the proof of the common-source theorem. First, let $X_{a,n_a}$ denote a sequence of random variables corresponding to the generation of datasets $\mathbf{X_{a,n_a}}$ where $n_a$ is the index that denotes the increasing number of subjects in the dataset with a fixed number of elements from within each subject.  Also, let $P_{\theta_a}^{n_a}$ denote the joint probability measure on $X_{a,n_a}$ for all $\theta_a \in \Theta_a$ including $\theta_{a_0}$.  The two sets of unknown-source observations, $\mathbf{x_b}$ and $\mathbf{x_c}$, will retain the previous notation since the size and the observation of these datasets will remain fixed. Now, let $\hat{\theta}_{a|M_2}$ denote the maximum likelihood estimator for $\theta_{a_0}$ that would be computed using the entire dataset  $\mathbf{X_{n_a}}$ when $\mathbf{x_b}$ and $\mathbf{x_c}$ are generated under the model $M_2$.  Finally, let $\beta_{cs}({X_{a,n_a}}, \mathbf{x_b}, \mathbf{x_c})$ denote the random function corresponding to the Bayes Factor given in Equation~\ref{Vcs_fact2} before the value of $X_{a,n_a}$ has been observed.

\begin{quote}\normalsize
\begin{thm*}[Common-Source Bayes Factor Consistency]
Given a fixed observation of $\mathbf{x_b}$ and $\mathbf{x_c}$, suppose that the likelihood ratio function $\lambda_{cs}(\theta_a|\mathbf{x_b}, \mathbf{x_c})$ is bounded in a neighborhood of $\theta_{a_0}$ and that $\hat{\theta}_{a|M_2}$ is a consistent estimator for $\theta_{a_0}$.  Furthermore, suppose the assumptions of the Bernstein-von Mises Theorem are satisfied.  Then as $n_a \to \infty$, the Bayes Factor converges in $P_{\theta_a}^{n_a}$-probability to the true likelihood ratio,
\begin{equation*}
\beta_{cs}({X_{a,n_a}}, \mathbf{x_b}, \mathbf{x_c})  \overset{P_{\theta_a}^{n_a}}{\longrightarrow}  \lambda_{cs}(\theta_{a_0}|\mathbf{x_b}, \mathbf{x_c}).
\end{equation*}
\end{thm*}
\end{quote}
\begin{proof}
First, let $X_{n_a}$ denote the random variable associated with generating the entire dataset where $X_{a,n_a}$ generates the observations from the population of sources and where the observations $\mathbf{x_b}$ and $\mathbf{x_c}$ are already fixed.  Also, let $\delta_{\theta_{a_0}}$ denote the probability measure that is degenerate at $\theta_{a_0}$. 
\begin{eqnarray*}
\Bigg|\beta_{cs}({X_{n_a}})  &-&  \lambda_{cs}(\theta_{a_0}|\mathbf{x_b}, \mathbf{x_c})\Bigg|
\\ &=& \Bigg|\int \lambda_{cs}(\theta_a| \mathbf{x_b}, \mathbf{x_c})\ d \left[ \Pi(\theta_a|{X_{n_a}}, M_2) - \delta_{\theta_{a_0}}(\theta_a) \right] \Bigg|
\end{eqnarray*}

Note that $\Pi(\theta_a|{X_{n_a}}, M_2) - \delta_{\theta_{a_0}}(\theta_a)$ is a sequence of signed measures.  For any signed measure $\mu$ it follows that $|\mu(A)| \leq ||\mu||_{TV}$ where $||\mu||_{TV}$ is the total variation norm \citet{Ash}.  Now, by the assumption that the likelihood ratio function is bounded, let $\lambda_{cs}(\theta_a| \mathbf{x_b}, \mathbf{x_c}) \leq C$ for some real number $C>0$. Therefore,
\begin{eqnarray*}
\Bigg|\beta_{cs}({X_{n_a}})  &-&  \lambda_{cs}(\theta_{a_0}|\mathbf{x_b}, \mathbf{x_c})\Bigg|
\\ &\leq& \int \lambda_{cs}(\theta_a| \mathbf{x_b}, \mathbf{x_c})\ \Bigg|d \left[ \Pi(\theta_a|{X_{n_a}}, M_2) - \delta_{\theta_{a_0}}(\theta_a) \right] \Bigg|
\\ &\leq& C\  \Bigg| \Bigg| \Pi(\theta_a|{X_{n_a}}, M_2) - \delta_{\theta_{a_0}}(\theta_a) \Bigg| \Bigg|_{TV}.
\end{eqnarray*}

Now, let $\Phi(\theta_a|{X_{n_a}}, M_2)$ denote the probability measure corresponding to the normal distribution with mean $\hat{\theta}_{a|M_2}$ and variance $\frac{1}{n_a}I_{\hat{\theta}_{a|M_2}}^{-1}$ where $I_{\hat{\theta}_{a|M_2}}^{-1}$ is the corresponding inverse of the observed Fisher's information matrix.  Then by the triangle inequality, we obtain
\begin{eqnarray*}
&\Bigg| \Bigg|& \hspace{-8pt} \Pi(\theta_a|{X_{n_a}}, M_2) - \delta_{\theta_{a_0}}(\theta_a) \Bigg| \Bigg|_{TV}
\\ &\leq& \Bigg| \Bigg| \Pi(\theta_a|{X_{n_a}}, M_2) -  \Phi(\theta_a|{X_{n_a}}, M_2)\Bigg| \Bigg|_{TV} + \Bigg| \Bigg|\Phi(\theta_a|{X_{n_a}}, M_2) - \delta_{\theta_{a_0}}(\theta_a) \Bigg| \Bigg|_{TV}.
\end{eqnarray*}
By the Bernstein-von Mises Theorem, as $n_a \to \infty$ then
\begin{equation*}
\Bigg| \Bigg| \Pi(\theta_a|{X_{n_a}}, M_2) -  \Phi(\theta_a|{X_{n_a}}, M_2)\Bigg| \Bigg|_{TV} \overset{P_{\theta_a}^{n_a}}{\longrightarrow} 0.
\end{equation*}
By the assumption that $\hat{\theta}_{a|M_2}$ is consistent and provided that $I_{\hat{\theta}_{a|M_2}}^{-1}$ is bounded in $P_{\theta_{a_0}}$-probability in a neighborhood of $\theta_{a_0}$, then this implies that as $n_a \to \infty$
\begin{equation*}
\Bigg| \Bigg|\Phi(\theta_a|{X_{n_a}}, M_2) - \delta_{\theta_{a_0}}(\theta_a) \Bigg| \Bigg|_{TV} \overset{P_{\theta_a}^{n_a}}{\longrightarrow} 0.
\end{equation*} 
\end{proof}

Please observe the following notational extensions to facilitate the proof of the specific-source theorem. First, let $X_{a,n_a}$ be defined in the same way as for the common-source model selection problem. Similarly, let $X_{b,n_b}$ denote a sequence of random variables corresponding to the generation of data from the fixed subject of interest where $n_b$ is the index that denotes the increasing number of elements from within that subject. For simplicity, we will fix $n_a = n_b$ so that the number of subjects in the population increases in the exact same way as the number of elements from the subject of interest.  The proofs of the results can easily be modified to accommodate more flexible relationships between the two sample sizes. In addition, let $P_{\theta}^n$ denote the joint probability measure on $X_{a,n_a}$ and $X_{b,n_b}$ for all $\theta \in \Theta$ where $\Theta$ is the joint parameter space, including $\theta_0$ which denotes the true value of the joint parameter.  The set of unknown-source observations will retain the notation $\mathbf{x_c}$ since the size and the observation of this set will remain fixed. Now, let $\hat{\theta}_{n|M_2}$ denote the maximum likelihood estimate computed using the entire dataset  $\mathbf{X_n}$ when $\mathbf{x_c}$ is generated under the model $M_2$.  Finally, let $\beta_{ss}({X_{a,n_a}}, X_{b,n_b}, \mathbf{x_c})$ denote the random function corresponding to the Bayes Factor given in Equation~\ref{Vss_fact2} before the values of $X_{a,n_a}$ and $X_{b,n_b}$ have been observed.

\begin{quote}\normalsize
\begin{thm*}[Specific-Source Bayes Factor Consistency]
Given a fixed observation of $\mathbf{x_c}$,  suppose that the likelihood ratio function, $\lambda_{ss}(\theta| \mathbf{x_c})$, is bounded in a neighborhood of $\theta_0$ and that $\hat{\theta}_{n|M_2}$ is a consistent estimator for $\theta_0$.  Furthermore, suppose the assumptions of the Bernstein-von Mises Theorem are satisfied. Then the Bayes Factor converges in $P_{\theta}^n$-probability to the true likelihood ratio as $n \to \infty$,
\begin{equation*}
\beta_{ss}({X_{a,n_a}}, X_{b,n_b}, \mathbf{x_c})  \overset{P_{\theta}^n}{\longrightarrow}  \lambda_{ss}(\theta_0| \mathbf{x_c}).
\end{equation*}
\end{thm*}
\end{quote}
\begin{proof}
Similar to the proof of Theorem~\ref{CS-BF-LR}.
\end{proof}

\section{Credible Interval Proofs}\label{C-ints}

Before we begin the proofs of the results for the credible intervals, it will be necessary to define some notational conventions related to the generation of the data.  To start, let $P_{\theta_a}^{n_a}$ denote the joint probability measure on $X_{a,n_a}$ for all $\theta_a \in \Theta_a$ including $\theta_{a_0}$ (for both the common-source and specific-source problems).  Similarly, for the specific-source problem only, let $P_{\theta_b}^{n_b}$ denote the joint probability measure on $X_{b,n_b}$ for all $\theta_b \in \Theta_b$ including $\theta_{b_0}$.  Next, let $X_{n}$ denote the random variable associated with generating the entire dataset under either of the non-nested model selection problems.  For the common-source problem, this means that $X_{a,n_a}$ generates the observations from the population of sources and where the observations $\mathbf{x_b}$ and $\mathbf{x_c}$ are already fixed.  For the specific source problem, this means that $X_{a,n_a}$ generates the observations from the larger population of subjects, $X_{b,n_b}$ generates the observations from the population associated with the single subject of interest, and where $\mathbf{x_c}$ is already fixed.  Furthermore, let $P_{\theta}^{n}$ denote the joint probability measure on $X_n$ for all $\theta \in \Theta$ including $\theta_{0}$ where $\Theta$ is the joint parameter space. 

Additional notational conventions needed to understand the proofs of the results for the credible intervals will now be considered. Let $\lambda(\theta)$ denote the likelihood ratio function for either the common-source or the specific-source problem given by Equation~\ref{LRcs} or Equation~\ref{LRss}, respectively.  Moreover, let $\mathcal{N}(\mu_n, \Sigma_n)$ denote the sequence of probability measures corresponding to the normal distribution with mean vector $\mu_n$ and a covariance matrix of $\Sigma_n$. Also, let $\hat{\theta}_{n|M_2}$ denote the maximum likelihood estimator computed from the entire collection of data $\mathbf{X_n}$ where the unknown-source set(s) of observations are generated according to the model $M_2$, and let $I_{\hat{\theta}_{n|M_2}}^{-1}$ denote the corresponding inverse of the observed Fisher's information matrix.  By properties of M-estimators, then $\hat{\theta}_{n|M_2}$ can be designed in a way such that it is a consistent estimator of $\theta_0$ \citet{VDV, VDVWellner}. Next, let, $\gamma_{n|M_2}^2$ be defined such that 
$$\gamma_{n|M_2}^2 = \lambda^{\prime}(\hat{\theta}_{n|M_2})^T  I_{\hat{\theta}_{n|M_2}}^{-1} \lambda^{\prime}(\hat{\theta}_{n|M_2})$$
where $\lambda^{\prime}$ is the vector of first partial derivatives of the likelihood ratio function. Finally, let $\beta({X_n})$ represent either the sequence of common-source or specific-source Bayes Factors.  It should be noted that further details of the proofs in this section are given in the supplementary material \citet{BFvsLRsupplement}.

\begin{quote}\normalsize
\begin{lem}\label{C-ints-1}
For a given observation of the unknown-source set(s) of observations, suppose that the likelihood ratio function, $\lambda(\theta)$, is twice continuously differentiable, and that $\hat{\theta}_{n|M_2}$ is a consistent estimator for $\theta_0$ under $P_{\theta}^n$-probability. If the assumptions of the Bernstein-von Mises Theorem hold, then $\Big|\Big| \Pi(\lambda(\theta)|{X_n}, M_2) -  \mathcal{N}(\lambda(\hat{\theta}_{n|M_2}), \frac{1}{n}\gamma_{n|M_2}^2) \Big|\Big|_{TV}$ converges to zero in $P_{\theta}^n$-probability as $n \to \infty$.
\end{lem}
\end{quote}
\begin{proof}
Consider the following result based on the Taylor series expansion of $\lambda(\theta)$ about the maximum likelihood estimate $\hat{\theta}_{n|M_2}$:
\begin{equation*}
\sqrt{n} \left[\lambda(\theta) - \lambda(\hat{\theta}_{n|M_2})\right] = \sqrt{n}(\theta-\hat{\theta}_{n|M_2})^T \lambda^{\prime}(\hat{\theta}_{n|M_2}) + \frac{\sqrt{n}}{2}(\theta-\hat{\theta}_{n|M_2})^T  \lambda^{\prime \prime}(\tilde{\theta}_{n|M_2}) (\theta-\hat{\theta}_{n|M_2})
\end{equation*}
where $\lambda^{\prime}$ is the vector of first partial derivatives of $\lambda$, $ \lambda^{\prime \prime}$ is the matrix of second partial derivatives of $\lambda$, and $\tilde{\theta}_{n|M_2}$ is a value on the line between $\theta$ and $\hat{\theta}_{n|M_2}$.  Now, consider the error term of the expansion given by
\begin{equation*}
\frac{1}{2} \left[ \sqrt{n} (\theta-\hat{\theta}_{n|M_2}) \right]^T \left[ \frac{1}{\sqrt{n}} \lambda^{\prime \prime}(\tilde{\theta}_{n|M_2}) \right] \left[ \sqrt{n}(\theta-\hat{\theta}_{n|M_2})\right].
\end{equation*}
The Bernstein-von Mises Theorem implies that the $\left[ \sqrt{n} (\theta-\hat{\theta}_{n|M_2}) \right]$ terms are bounded in $P_{\theta}^n$-probability as $n \to \infty$.  Under the assumption that the maximum likelihood estimate $\hat{\theta}_{n|M_2}$ is consistent, i.e. $\hat{\theta}_{n|M_2} \overset{P_{\theta_{0}}}{\longrightarrow} \theta_0$ as $n \to \infty$, then the Continuous Mapping Theorem and the Squeeze Theorem imply that the $\left[ \frac{1}{\sqrt{n}} \lambda^{\prime \prime}(\tilde{\theta}_{n|M_2}) \right]$ term converges in $P_{\theta}^n$-probability to zero as $n \to \infty$.  Therefore, the error term converges in $P_{\theta}^n$-probability to zero  as $n \to \infty$.

This implies that the distribution of $\sqrt{n} \left[\lambda(\theta) - \lambda(\hat{\theta}_{n|M_2})\right]$ is determined by the the remaining term in the expansion, $\sqrt{n}(\theta-\hat{\theta}_{n|M_2})^T \lambda^{\prime}(\hat{\theta}_{n|M_2})$.  By the Bernstein-von Mises Theorem, then as $n \to \infty$
\begin{equation*}
\Big|\Big| \Pi(\sqrt{n}(\theta-\hat{\theta}_{n|M_2})|{X_n}, M_2) - \mathcal{N}(0, I_{\hat{\theta}_{n|M_2}}^{-1}) \Big|\Big|_{TV} \overset{P_{\theta}^n}{\longrightarrow} 0.
\end{equation*}
Now, by properties of normal distributions regarding linear combinations of normal random variables, then this implies that as $n \to \infty$
\begin{equation*}
\Big|\Big| \Pi(\sqrt{n}(\theta-\hat{\theta}_{n|M_2})^T \lambda^{\prime}(\hat{\theta}_{n|M_2})|{X_n}, M_2) - \mathcal{N}(0, \gamma_{n|M_2}^2) \Big|\Big|_{TV} \overset{P_{\theta}^n}{\longrightarrow} 0.
\end{equation*}
Lastly, Slutsky's Lemma and the re-centered and rescaled version of the Bernstein-von Mises Theorem result gives, as $n \to \infty$,
\begin{equation*}
\Big|\Big| \Pi(\lambda(\theta)|{X_n}, M_2) - \mathcal{N}(\lambda(\hat{\theta}_{n|M_2}), \frac{1}{n}\gamma_{n|M_2}^2) \Big|\Big|_{TV} \overset{P_{\theta}^n}{\longrightarrow} 0.
\end{equation*}
\end{proof}

\begin{quote}\normalsize
\begin{lem}\label{C-ints-2}
For a given observation of the unknown-source set(s) of observations, suppose that $\hat{\theta}_{n|M_2}$ is a consistent estimator for $\theta_0$ under $P_{\theta}^n$-probability and that $\sigma_n^2 \to 0$ as $n \to \infty$.  If the assumptions of the Bernstein-von Mises Theorem hold, then $\Big|\Big|\mathcal{N}(\lambda(\hat{\theta}_{n|M_2}), \frac{1}{n}\gamma_{n|M_2}^2) - \mathcal{N}(\beta({X_n}), \sigma_n^2) \Big|\Big|_{TV}$ converges to zero in $P_{\theta}^n$-probability as $n \to \infty$.
\end{lem}
\end{quote}
\begin{proof}
First, we will show that the difference in the means of the two sequences of Normal distributions will converge to 0. By Theorem~\ref{CS-BF-LR} and Theorem~\ref{SS-BF-LR}, we have that $\beta({X_n})$ converges in $P_{\theta}^n$-probability to the true likelihood ratio, $\lambda(\theta_0)$, as $n \to \infty$ for both the common-source and the specific-source Bayes Factors.  Now as $n \to \infty$, it readily follows that, for both the common-source and the specific-source likelihood ratio function,
\begin{equation*}
\lambda(\hat{\theta}_{n|M_2}) \overset{P_{\theta}^n}{\longrightarrow} \lambda(\theta_0).
\end{equation*}
Secondly, we will show that the difference between the variances of the two sequences of Normal distributions will converge to 0.  By assumption, $\sigma_n^2 \to 0$ as $n \to \infty$. Next, under the assumptions of the Bernstein-von Mises Theorem, $\frac{1}{n}\gamma_{n|M_2}^2$ also converges to 0 as $n \to \infty$.
Finally, the Continuous Mapping Theorem gives the needed result:
\begin{equation*}
\Big|\Big| \mathcal{N}(\lambda(\hat{\theta}_{n|M_2}), \frac{1}{n}\gamma_{n|M_2}^2) - \mathcal{N}(\beta({X_n}), \sigma_n^2) \Big|\Big|_{TV} \overset{P_{\theta}^n}{\longrightarrow} 0.
\end{equation*}

\end{proof}

\begin{quote}\normalsize
\begin{thm*}[Approximate $1-\alpha$ Credible Interval for the LR]
Let let the assumptions of Lemma~\ref{C-ints-1} and Lemma~\ref{C-ints-2} hold and let 
$$\mathcal{I}_n = \beta({X_n}) \pm \Phi^{-1}(1-\alpha/2)\sigma_n$$
where $\beta({X_n})$ represents the sequence of either the common-source or the specific-source Bayes Factor, $0<\alpha<1$ is the desired significance level, $\sigma_n$ is the sequence of posterior standard errors for the likelihood ratio, and $\Phi^{-1}$ is the standard normal quantile function. Then as $n \to \infty$
\begin{equation*}
	\Pi(\lambda(\theta) \in \mathcal{I}_n |{X_n}, M_2) \to 1-\alpha.
\end{equation*}

\end{thm*}	
\end{quote}
\begin{proof}

By Lemma~\ref{C-ints-1} and Lemma~\ref{C-ints-2}, we have that as $n \to \infty$
\begin{equation}\label{step3}
\Big|\Big| \Pi(\lambda(\theta)|{X_n}, M_2) -  \mathcal{N}(\beta({X_n}), \sigma_n^2) \Big|\Big|_{TV} \overset{P_{\theta}^n}{\longrightarrow} 0.
\end{equation}
Since $\mathcal{N}(\beta({X_n}), \sigma_n^2)$ denotes the sequence of probability measures corresponding to the normal distribution with mean $\beta({X_n})$ and covariance $\sigma_n^2$, then
$$\mathcal{I}_n = \beta({X_n}) \pm \Phi^{-1}(1-\alpha/2)\sigma_n$$
is a sequence of corresponding credible intervals with $100(1-\alpha)$\% credibility. The previous result in Equation~\ref{step3} then implies that as $n \to \infty$
\begin{equation*}
\Pi(\lambda(\theta) \in \mathcal{I}_n|{X_n}, M_2)  \to 1-\alpha.
\end{equation*}

\end{proof}

\end{document}